\documentclass[a4paper,10pt]{amsart}
\usepackage{verbatim,latexsym,exscale,amssymb,amsthm,amsmath,amsfonts,mathrsfs,amsrefs}
\usepackage{stmaryrd}
\usepackage[all]{xy}

\DeclareMathSymbol{\Gamma}{\mathalpha}{letters}{"00}

\xyoption{ps}
\SelectTips{cm}{}

\newdir{ >}{{}*!/-6pt/@{>}}

\newdir{ >}{{}*!/-5pt/\dir{>}}

\newcommand{\tcat}{\mathcal{T}}

\newcommand{\spec}{\mathop{\mathrm{Spec}}}

\newcommand{\Spec}{\mathop{\mathrm{Spec}}}

\newcommand{\supp}{\mathrm{supp}}
\newcommand{\Supp}{\mathop{\mathrm{Supp}}}

\newcommand{\Hom}{\mathrm{Hom}}

\newcommand{\End}{\mathrm{End}}

\newcommand{\N}{\mathbb{N}}

\newcommand{\cone}{\mathop{\mathrm{cone}}}

\newcommand{\colim}{\mathop{\mathrm{colim}}}
\newcommand{\hocolim}{\mathop{\mathrm{hocolim}}}
\newcommand{\unit}{\mathbf{1}}

\newcommand{\Ker}{\mathop{\mathsf{Ker}}}

\newcommand{\Cstar}{\mathrm{C}^*\!}

\newcommand{\Ext}{\mathrm{Ext}}
\newcommand{\Tor}{\mathrm{Tor}}

\newcommand{\Kth}{K}
\newcommand{\KK}{\mathsf{KK}}

\newcommand{\Boot}{\mathsf{Boot}}

\newcommand{\D}{\mathsf{D}}

\newcommand{\Ab}{\mathsf{Ab}}

\newcommand{\Z}{\mathbb{Z}}
\newcommand{\R}{\mathbb{R}}
\newcommand{\C}{\mathbb{C}}
\newcommand{\Q}{\mathbb{Q}}

\theoremstyle{definition}

\newtheorem{defi}{Definition}[section]

\newtheorem*{conv*}{Conventions}
\newtheorem*{ack}{Acknowledgement}

\theoremstyle{theorem}

\newtheorem{thm}[defi]{Theorem}
\newtheorem{lemma}[defi]{Lemma}
\newtheorem{thm-defi}[defi]{Theorem-Definition}
\newtheorem{prop}[defi]{Proposition}
\newtheorem{cor}[defi]{Corollary}

\newtheorem*{thm*}{Theorem}
\newtheorem*{lemma*}{Lemma}
\newtheorem*{cor*}{Corollary}
\newtheorem*{conj*}{Conjecture}

\theoremstyle{remark}

\newtheorem{notation}[defi]{Notation}
\newtheorem{remark}[defi]{Remark}
\newtheorem*{remark*}{Remark}

\begin{document}
\title[Localizing subcategories of the Bootstrap category]{Localizing subcategories in the Bootstrap category of separable $\Cstar$-algebras}
\author{Ivo Dell'Ambrogio}

\thanks{Research supported by the Stefano Franscini Fund, Swiss National Science Foundation grant Nr.\ PBEZP2-125724.}
\address{Dpt.~of Mathematics, Faculty of Science, National University of Singapore (NUS), 10~Lower Kent Ridge Road, S17--08--22, Singapore, 117543.}
\email{ivo@dellambrogio.ch}

 \date{}

\begin{abstract}
Using the classical universal coefficient theorem of Rosenberg-Schochet, we prove a simple classification of all localizing subcategories of the Bootstrap category $\Boot\subset \KK$ of separable complex C*-algebras. Namely, they are in a bijective correspondence with subsets of the Zariski spectrum $\spec \Z$ of the integers -- precisely as for the localizing subcategories of the derived category $\D(\Z)$ of complexes of abelian groups.
We provide corollaries of this fact and put it in context with the similar classifications available in the literature. 
\end{abstract}

\maketitle


\section{Introduction}

Denote by $\KK$ the Kasparov category of separable complex C*-algebras, with Hom groups given by G.\ Kasparov's bivariant K-groups $\KK(A,B)$ and with composition the Kasparov product (\cite{kasparov_first,blackadar}).
As observed by R.\ Meyer and R.\ Nest \cite{meyernest_bc}, $\KK$ carries the structure of a triangulated category (see \cite{neemanTr}), which conveniently captures the homological properties of KK-theory in a well-understood abstract framework. The Kasparov category also has arbitrary countable coproducts, given by the C*-algebra direct sums. The \emph{Bootstrap category}, which we denote by $\Boot$, can then be defined as the smallest triangulated subcategory\footnote{All subcategories considered here will be full and closed under isomorphic objects.} of $\KK$ which contains $\C$ and is closed under the formation of coproducts and retracts of objects. More generally, a subcategory $\mathcal L$ of a triangulated category $\mathcal T$ is said to be \emph{localizing} if it is a triangulated subcategory closed under the formation of retracts and whatever coproducts exist in $\tcat$ (if $\tcat$ has at least all countable coproducts, as here, the closure under taking retracts follows then from the other properties).  

The interest in the Bootstrap category is due to the classical result of J.\ Rosenberg and C.\ Schochet \cite{rs}, that the C*-algebras in it satisfy a universal coefficient theorem (UCT), recalled below together with its immediate corollaries. The UCT shows that $\Boot$ is not far from the category of $\Z/2$-graded countable abelian groups. In this short note we use the UCT, and the related K\"unneth theorem, to prove a very simple classification of all localizing subcategories of the Bootstrap category, as follows.
For a C*-algebra $A$, denote by
$
K_*(A; \mathbb F_p)
$
 its topological K-theory with coefficients in the residue field $\mathbb F_p$ at $p \in \spec \Z$. That is, $\mathbb F_p= \Z/p$ if $p$ is a prime number and $\mathbb F_p=\Q$ for $p=0$; recall that we may set $K_*(A;\mathbb F_p):= K_*(A\otimes \kappa(p))$, where $\kappa(p)$ is any C*-algebra in $\Boot$ with $K_0(\kappa(p))=\mathbb F_p$ and $K_1(\kappa(p))=0$. 
 
 Here is the classification:

\begin{thm} \label{thm:main}
There is an inclusion-preserving bijection between localizing subcategories of the Bootstrap category $\Boot$ on the one hand, and subsets of the Zariski spectrum $\Spec \Z$ on the other, given by the maps
\begin{eqnarray*}
\mathcal Loc(\Boot) &\simeq & \mathcal P( \spec \Z ) \\
\mathcal L &\mapsto  & \Supp \mathcal L:= \{ p \mid \exists\, A\in \mathcal L \textrm{ s.\,t.\ }  K_*(A; \mathbb F_p) \neq0 \} \\
\{A\mid  \forall p \not\in S, K_*(A; \mathbb F_p)= 0 \} =: \mathcal L_S &\mapsfrom & S
\end{eqnarray*} 
 This bijection restricts to smashing subcategories of $\Boot$ on one side and specialization closed subsets of $\spec \Z$ on the other side.
\end{thm}

In this note, a localizing subcategory $\mathcal L$ is  \emph{smashing} if its right orthogonal
\begin{equation*}
\mathcal L^{\perp} := \{ A\in \Boot \mid \KK(B,A)=0 \textrm{ for all } B\in \mathcal L \}
\end{equation*} 
is again a localizing subcategory of $\Boot$. Thus, for instance, the subcategory $\mathcal L_{\{0\}}$ of C*-algebras whose mod~$p$ K-theory vanishes for all prime numbers~$p$, is localizing but not smashing. Note also that a subset $S\subseteq \spec \Z$ is specialization closed if and only if it is either a set of prime numbers or the whole spectrum.

Theorem \ref{thm:main} and its proof have some easy corollaries. For any family $\mathcal E \subset \Boot$ of object, let $\langle \mathcal E\rangle $ denote the localizing subcategory  \emph{generated by $\mathcal E$}, i.e., the smallest localizing subcategory of $\Boot$ containing $\mathcal E$. Thus for example $\Boot=\langle \C\rangle$. Also, the localizing subcategory corresponding in Thm.~\ref{thm:main} to the set $S\subseteq \spec \Z$ can be described as $\mathcal L_S=\langle \kappa(p)\mid p\in S \rangle$ (see the beginning of~$\S$\ref{subsec:loc_subcats}).

An object $C\in \Boot$ is \emph{compact} (more precisely we should say \emph{compact${}_{\aleph_1}$}, see \cite{kkGarticle}*{\S 2.1}) if its Hom functor $\KK(C,-)$ takes values in countable abelian groups and commutes with coproducts.
Write $\Boot_c$ for the subcategory of compact objects in $\Boot$; it is a thick triangulated subcategory, consisting precisely of all algebras $A\in \Boot$ with finitely generated K-theory groups (Lemma \ref{lemma:cptobj}).

The next corollary is the analog for $\Boot$ of the telescope conjecture in stable homotopy theory, or in more general compactly generated categories (see e.g.~\cite{krause:telescope}).

\begin{cor}[Telescope conjecture for $\Boot$] \label{cor:telescope}
Every smashing subcategory $\mathcal L$ of $\Boot$ is generated by the compact objects it contains: $\mathcal L=\langle \mathcal L\cap \Boot_c \rangle\subseteq \Boot$.
\end{cor}

\proof
Clearly $\langle \mathcal L\cap \Boot_c\rangle \subseteq \mathcal L$, since $\mathcal L$ is localizing. If $\mathcal L$ is moreover smashing, then by our classification it is either $\mathcal L_{\spec \Z} = \Boot$, which is generated by the compact object~$\C$, or it is $\mathcal L_S$ for some set $S$ of prime numbers, and in this case $\mathcal L_S$ is generated by the set $\{\kappa(p)\mid p\in S\} \subseteq \mathcal L\cap \Boot_c$ of its compact objects.  
\qed

\begin{cor} \label{cor:thick}
There is a bijection between thick subcategories of $\Boot_c$ (i.e., triangulated subcategories closed under taking retracts of objects) and specialization closed subsets of $\Spec \Z$.
\end{cor}

\proof
This is immediate from Cor.~\ref{cor:telescope} and the last statement of Thm.~\ref{thm:main}, once we know that the maps $\mathcal C\mapsto \langle \mathcal C\rangle$ and $\mathcal L\mapsto \mathcal L\cap \Boot_c$ provide a bijection between thick subcategories $\mathcal C\subseteq \Boot_c$ of compact${}_{\aleph_1}$ objects and localizing subcategories $\mathcal L\subseteq \Boot$ which are generated by compact${}_{\aleph_1}$ objects (because $\Boot$ is compactly${}_{\aleph_1}$ generated, see \cite{kkGarticle}*{$\S2.1$}).
\qed

If we also take into consideration the biexact tensor product $\otimes=\otimes_{\mathrm{min}}$, we can easily derive from the bijection of Cor.~\ref{cor:thick} a new proof of \cite{kkGarticle}*{Thm.~6.9}:

\begin{cor} \label{cor:spec_bootc}
There is a canonical isomorphism of schemes
$$\spec(\Boot_c)\simeq \spec(\Z),$$
where the first $\Spec$ denotes P.\ Balmer's spectrum of prime $\otimes$-ideals \cite{balmer_prime}. 
\qed
\end{cor}

%

We leave this to the interested reader.

In the remaining pages, we shall prove Theorem \ref{thm:main} using general results on triangulated categories, together with Bott periodicity and the UCT. The next challenge now would be to provide a similar classification for the strictly bigger (\cite{skandalis_nucl}*{p.\ 751}) triangulated category $\KK$, where the UCT does not hold. Indeed, not much is known about $\KK$ -- structurally speaking -- except that it ``decomposes'' into the orthogonal localizing subcategories $\Boot$ and $\Ker K_* := \{A \in \KK \mid K_*A=0 \}=\Boot^{\perp}$; that is, every $A\in \KK$ fits into an (automatically unique and functorial) exact triangle $A'\to A\to A'' \to \Sigma A$ with $A'\in \Boot$ and $A''\in \Ker K_*$. This is due to Meyer and Nest (\cite{meyernest_hom,meyer_hom}; cf.\ \cite{kkGarticle}*{Thm.~5.8}). Hence structural questions about $\KK$ could hopefully be reduced to questions about $\Boot$, which is well-understood, and about $\Ker K_*$, which is still mysterious. At this point, one is led to investigate homological invariants of separable C*-algebras finer than K-theory.
 
In another direction, we hope to prove similar results for a suitable `bootstrap subcategory' of the $G$-equivariant Kasparov category for a finite group~$G$, where the role of the integers is taken over by the complex character ring of~$G$. 
 
\section{The universal coefficient and K\"unneth theorems}

Let us recall some well-known results, mostly for the purpose of fixing notation.

As is customary, we write $\KK_*(A,B)$ for the $\Z/2$-graded abelian group with $\KK(A,B)$ in degree zero and $\KK(\Sigma A,B)$ in degree one. Here $\Sigma= C_0(\R)\otimes- $ is the suspension functor of the triangulated structure of $\KK$. The Bott isomorphism $C_0(\R)\otimes C_0(\R)\simeq \C$ in $\KK$ shows that $\Sigma\simeq \Sigma^{-1}$. The K-theory functor $K_*=\KK_*(\C,-)$ is a coproduct preserving homological functor taking values in the abelian category $\Ab^{\Z/2}_{\infty}$ of $\Z/2$-graded countable abelian groups and degree-preserving homomorphisms, for which we shall use the following 

\begin{notation} \label{notation:graded}
If $M$ is a $\Z/2$-graded abelian group, then $M_0$ shall denote its zero-degree part and $M_1$ its one-degree part; 
on the other hand, if $M$ is an abelian group and $\varepsilon\in \Z/2$, we write $M[\varepsilon]$ for the graded group with $M$ in degree $\varepsilon$ and $0$ in degree~$\varepsilon+1$.
\end{notation}

\subsection{The universal coefficient theorem and some corollaries}
The extremely useful UCT of Rosenberg-Schochet \cite{rs} computes the Kasparov groups $\KK_*(A,B)$ in terms of the K-theory groups of $A$ and $B$:

\begin{thm}[Universal coefficient theorem] \label{thm:UCT}
Let $A$ and $B$ be separable $\Cstar$-al\-ge\-bras, with $A\in\Boot$. Then there is a short exact sequence of $\Z/2$-graded groups
\begin{equation} \label{uct}
\xymatrix{ \Ext(\Kth_*A, \Kth_*B)\ar@{  >->}[r]^-{+1} & \KK_*(A,B)\ar@{->>}[r] & \Hom(\Kth_*A, \Kth_*B),}
\end{equation}
where the map marked by $+1$ has degree one, and the second map is the $\Kth$-theory functor $\Kth_0=\KK(\C,-)$.
The sequence is natural and splits unnaturally. 
\qed
\end{thm}

Here $\Hom$ and $\Ext$ denote the (graded) Hom and Ext groups computed in $\Ab^{\Z/2}_{\infty}$. 
Concretely, $\Hom$ has the degree preserving homomorphisms in degree zero and the degree exchanging homomorphisms in degree one, while 
$\Ext(K_*A,K_*B)_0= \Ext^1_{\Z}(K_0A,K_0B)\oplus \Ext^1_{\Z}(K_1A,K_1B)$ and 
$\Ext(K_*A,K_*B)_1= \Ext^1_{\Z}(K_0A,K_1B)\oplus \Ext^1_{\Z}(K_1A,K_0B)$
in terms of the usual $\Ext^1_{\Z}$ of abelian groups. 

\begin{remark}
The UCT, and the K\"unneth theorem recalled below, are best understood in the general context of relative homological algebra for triangulated categories. There they are corollaries, among others, of the stronger result that homological algebra in $\KK$ relative to the (homological ideal defined by) K-theory is \emph{hereditary}. Consult  \cite{meyernest_hom, meyer_hom} for this and more.  
\end{remark}

\begin{remark}
In fact, the UCT has the following converse: for any separable C*-algebra $A$, the sequence \eqref{uct} is exact (for all $B\in \KK$) if and only if $A\in \Boot$. This can be used to show that all \emph{commutative} separable C*-algebras belong to $\Boot$ (\cite{skandalis_nucl}*{Prop.\ 5.3}). The same question for \emph{nuclear} separable C*-algebras, however, is still unanswered. Indeed, this has been called ``one of the outstanding open questions of C*-algebra theory'' (\cite{blackadar}*{23.15.12}; see also \cite{dadarlat}).
\end{remark}

The following corollaries are easily derivable from the UCT and the 6-periodic long exact sequence in K-theory, see \cite{blackadar}*{\S 23} or \cite{kkGarticle}*{Cor.~5.29-32}.

\begin{cor} \label{cor:iso_proj_inj}
K-theory induces an isomorphism $\KK(A,B)\simeq \Hom(K_*A,K_*B)_0$ for any $A\in \Boot$ and $B\in \KK$ such that $K_*A$ is free or $K_*B$ is divisible.
\qed
\end{cor}

\begin{cor} \label{cor:realizeK}
Let $M$ be any $\Z/2$-graded countable abelian group. Then there is a C*-algebra $A\in\Boot$ in the Bootstrap category with $\Kth_*A\simeq M$.
\qed
\end{cor}

With this corollary, for instance, one constructs the objects $\kappa(p)$ with $K_*(\kappa(p))=\mathbb F_p[0]$ used in the definition of K-theory with coefficients. For $p\neq 0$, of course, one needs simply set $\kappa(p)= \cone(p)$, where $p\in \Z=\End_{\Boot}(\C)$ is considered to be a map in $\Boot$ and $\cone(p)$ is its cone (or equivalently: $\cone(p)=\Sigma C_p$, where $C_p$ is the ``mapping cone'', in the traditional arrow-reversing sense of C*-algebraists, of the unital *-homomorphism $\C\to M_p(\C)$).

%

\begin{cor} \label{cor:isotype}
Consider objects $A,B\in\Boot$ and an isomorphism $f:K_*A\simeq K_*B$ of graded abelian groups. Then $f$ lifts to an isomorphism $A\simeq B$ in $\Boot$.
\qed
\end{cor}

Thus K-theory induces a bijection between the isomorphism types of $\Boot$ and those of $\Ab^{\Z/2}_{\infty}$. It is no surprise that something like Theorem \ref{thm:main} must be true. 

\begin{remark} \label{rem:aleph_1}
As already mentioned in the Introduction, $\Boot$ is \emph{compactly${}_{\aleph_1}$ generated} in the sense of \cite{kkGarticle}, that is, it admits all countable coproducts and it is generated by a set of \emph{compact${}_{\aleph_1}$} objects. Note that $\Boot$ is not compactly generated in the usual sense, because it has no non-trivial uncountable coproducts -- this is due to the separability hypothesis built into the theory. This said, the reader may now relax in view of the following characterization.
\end{remark}

%

%

\begin{lemma}\label{lemma:cptobj}
An object $A\in \Boot$ is compact (i.e., compact${}_{\aleph_1}$) if and only if $\Kth_*A$ is finitely generated.
\end{lemma}

\proof 
One implication, that every compact object in $\Boot$ has finitely generated K-theory, follows from a routine induction based on the fact that the subcategory of compact objects coincides with the thick triangulated subcategory of $\Boot$ generated by $\C$.
If $A\in \Boot$ has finitely generated K-theory and $\coprod_iB_i$ is a coproduct in $\Boot$ (or even in $\KK$), then the canonical map $\bigoplus_i \KK(A,B_i)\to \KK(A,\coprod_iB_i)$ is an isomorphism by \cite{rs}*{Prop.\ 7.13} (which is another easy corollary of the UCT).
Again by the UCT, we see that $\KK(A,B)$ is countable whenever $K_*(A)$ is finitely generated (we also see that it is \emph{not} countable in general, e.g.\ $B=\C$ and $A=\coprod_{\N}\C$). 
\qed

\subsection{The K\"unneth theorem and some corollaries}
The Kasparov category $\KK$ comes equipped with a tensor product (a symmetric monoidal structure), induced by the minimal tensor product $\otimes:=\otimes_{\min}$ of C*-algebras, which preserves coproducts and exact triangles in each variable. 
The following K\"unneth theorem (not to be confused with the homonymous and allied, but different, K\"unnneth theorem \cite{rs}*{1.18}) computes the K-theory of a tensor product in terms of that of its factors, provided one of them is in $\Boot$.

\begin{thm}[K\"unneth theorem; \cite{blackadar}*{\S 23}, \cite{meyernest_hom}] \label{thm:KTP}
Let  $A$ and $B$ be separable $\Cstar$-al\-ge\-bras with $A\in\Boot$. Then there is a natural short exact sequence of $\Z/2$-graded abelian groups 
\begin{equation} \label{ktp}
\xymatrix{ \Kth_*(A)\,\hat{\otimes}\,\Kth_*(B)\ar@{  >->}[r]&\Kth_*(A\otimes B)\ar@{->>}[r]^-{+1}&\Tor(\Kth_*(A),\Kth_*(B))}
\end{equation}
which splits unnaturally. The first map is the one induced by the symmetric monoidal structure $\otimes:=\otimes_{\min}:\KK(\C,A)\times \KK(\C,B)\to\KK(\C,A\otimes B)$ of $\KK$.
\qed
\end{thm}

Here $\hat\otimes$ is the tensor product of $\Z/2$-graded groups, defined by 
$(M\,\hat\otimes\, N)_{\varepsilon} = \bigoplus_{i+j=\varepsilon} M_i \oplus M_j$, and $\Tor$ is its first derived functor, or concretely
$\Tor(M,N)_{\varepsilon}= \bigoplus_{i+j=\varepsilon} \Tor^{\Z}_1(M_i,N_j)$,
in terms of the usual $\Tor^{\Z}_1$ of abelian groups ($i,j,\varepsilon \in \Z/2$).

\begin{cor} \label{cor:tensor_dichotomy}
For all $p,q\in\spec \Z$, we have $\kappa(p)\otimes \kappa(q)\neq 0$ if and only if $p=q$.
\end{cor}

\proof
Immediate from the K\"unneth theorem, since $\Tor_*^{\Z}(\mathbb F_p,\mathbb F_q)\neq 0$ iff $p=q$.
\qed

\begin{cor} \label{cor:supp}
If $A\in \langle \kappa(p)\mid p\in S\rangle$ for some $S\subseteq \spec \Z$, then $A\otimes \kappa(q)=0$ for all $q\in \spec\Z \smallsetminus S$.
\end{cor}

\proof
Because the kernel on objects $\Ker(-\otimes \kappa(q))$ of tensorization with $\kappa(q)$ is a localizing subcategory of $\Boot$ and contains $\kappa(p)$ for $p\neq q$ by Cor.~\ref{cor:tensor_dichotomy}. 
\qed

\begin{cor} \label{cor:residual_decomp}
For every $A\in \Boot$ and $p\in \spec \Z$, the tensor product $A\otimes \kappa(p)$ is isomorphic to a coproduct of suspensions of copies of $\kappa(p)$. 
\end{cor}

\proof
By the K\"unneth theorem, $K_*(A\otimes \kappa(p))$ is a graded $k(p)$-vector space. Hence $A\otimes \kappa(p)$ has the desired form because of Cor.~\ref{cor:isotype}.
\qed

\subsection{Injective objects} \label{subsec:inj}
We recall a few well-known facts. The ring $\Z$ of integers is a commutative noetherian ring of homological dimension one, and it follows that every abelian group $M$ has a length one \emph{minimal} injective resolution (one which is a retract of any other). Such a resolution $M \rightarrowtail I^0 \twoheadrightarrow I^1$ is unique up to isomorphism.
The class of injective $\Z$-modules coincides with that of divisible abelian groups, and every divisible group $M$ is isomorphic to a direct sum of copies of the following groups -- the indecomposable injective $\Z$-modules:
\begin{equation*}
I(0):= \Q  \quad \textrm{ and } \quad I(p):= \Z[\textstyle \frac{1}{p}]/\Z   \quad, \quad  \textrm{ for } p \textrm{ prime } \text .
\end{equation*}
(This description is a simple special case of the well-known classification of injective modules over a commutative noetherian ring~\cite{matlis}.) Note that $I(p)$ consists entirely of $p$-primary torsion elements ($p\neq 0$).
If we localize the $\Z$-module $I(p)$ at $q$, we get $I(p)_{(q)}=0$ except if $q=p$ or $p=0$, in which case $I(p)_{(q)}=I(p)$.

All of this extends to graded abelian groups:
it is clear that an object $M \in \Ab_{\infty}^{\Z/2}$ is injective iff its zero- and one-degree parts are both (countable) injective abelian groups, and therefore $M\simeq \bigoplus_{n\in N_0}I(p_n)[0]\oplus \bigoplus_{m\in N_1}I(p_m)[1]$
for some countable sets $N_0$ and $N_1$. 

If we also factor in Corollary \ref{cor:iso_proj_inj}, we get
\begin{cor} \label{cor:iso_inj}
The K-theory functor $K_*:\Boot \to \Ab^{\Z/2}_{\infty}$ restricts to an equivalence between  the full subcategory of those $A\in \Boot$ with divisible K-theory, and the full subcategory of injective objects in $\Ab^{\Z/2}_{\infty}$.
\qed
\end{cor}

\begin{defi}
As we did for the \emph{residue field objects} $\kappa(p)$, we can use Cor.~\ref{cor:realizeK} to define for each $p\in \spec \Z$ the -- up the isomorphism, unique -- object  $\iota(p)$  in $\Boot$  such that
\begin{equation*}
K_*(\iota(p)) \simeq I(p)[0],
\end{equation*}
thus realizing in the Bootstrap category the indecomposable injectives $I(p)[0]$ and $I(p)[1]$ of $\Ab^{\Z/2}_{\infty}$ associated to~$p$. 
Note that $\iota(0)=\kappa(0)$ (recall that by definition $\mathbb F_0=\Q$ and thus $K_*(\kappa(0))=\Q[0]$), and that none of the $\iota(p)$ is compact.

\end{defi}

\begin{cor} \label{cor:inj_decomp}
Every $B\in \Boot$ with divisible $K$-theory is a coproduct of copies of $\iota(p)$ and $\Sigma\iota(p)$, $p\in \spec \Z$. 
\qed
\end{cor}

\subsection{Re-generating the Bootstrap category} We now investigate more closely the relation between the objects $\kappa(p)$ and $\iota(p)$. As it turns out, they generate the same localizing subcategories of $\Boot$ (Prop.\ \ref{prop:generation}).

\begin{defi}[Cf.\ \cite{bik}]  \label{defi:BootV}
Define the \emph{support}, written $\supp_{\Z} M$, of a (graded) abelian group as the set of $p\in \spec \Z$ such that $I(p)$ appears (in some degree) as a direct summand in a minimal injective resolution $I^0\to I^1$ of $M$.
For any subset $S\subseteq \spec \Z$, define 
\begin{equation*}
\Boot_S:= \{A\in \Boot \mid \supp_{\Z} K_*A \; \subseteq \; S \}.
\end{equation*}
\end{defi}

Thus for instance $\Boot_{\emptyset}=0$, $\Boot_{\spec \Z}=\Boot$, and $\Boot_{\spec \Z\smallsetminus \{0\}}$ consists of all C*-algebras with torsion K-theory. More generally:

\begin{lemma} \label{lemma:localizing}
If $V\subseteq \spec \Z$ is specialization closed (i.e., either $V=\spec \Z$ or $V$ is a subset of prime numbers), then 
$\Boot_V= \Ker \big(\bigoplus_{q\not\in V} K_*(-)_{(q)} \big) \subseteq \Boot$. In particular, as the object-kernel of a coproduct preserving homological functor, $\Boot_V$ is a localizing subcategory of $\Boot$.
\end{lemma}

\proof
This is \cite{bik}*{Lemma 2.3 (1)} for $R=\Z$. Of course, this case is also easy to see directly from $\S\ref{subsec:inj}$.
\qed

\begin{remark}
The hypothesis that $V$ be specialization closed is necessary. Indeed, if $S\subseteq \spec \Z$ is not specialization closed then $S$ contains $0$ without containing some prime number~$p$. Therefore $\kappa(0)\in \Boot_{\{0\}}\subseteq \Boot_S$ but $(K_* \kappa(0))_{(p)}=\Q[0]\neq 0$.
\end{remark}

\begin{lemma} \label{lemma:easy_generation}
$\Boot_V= \langle \iota(p)\mid p\in V \rangle$ for any specialization closed $V\subseteq \spec \Z$.
\end{lemma}

\proof 
By definition, for every $A\in \Boot_V$ there is in $\Ab^{\Z/2}_{\infty}$ a short exact sequence $K_*A \rightarrowtail I^0\twoheadrightarrow I^1$ such that $I^0$ and $I^1$ are direct sums of copies of $I(p)[0]$ and $I(p)[1]$, for $p\in V$. 
Now use the UCT to realize $f:I^0\to I^1$ as a map $\varphi:B^0\to B^1$ in $\Boot$. We see from $\S$\ref{subsec:inj} that $B^0$ and $B^1$ must be isomorphic to coproducts of suspensions of copies of the objects $\iota(p)$, $p\in V$. In particular $B^0,B^1 \in \langle \iota(p)\mid p\in V \rangle\subseteq \Boot$. Finally, the 6-periodic long exact sequence in K-theory for a distinguished triangle $\Sigma B^1 \to C \to B^0 \to B^1$ containing $\varphi$ splits in two short exact sequences by the surjectivity of $K_*(\varphi)=f$, and we immediately see that $K_*A\simeq K_*C$ and therefore $A\simeq C\in \langle \iota(p)\mid p\in V \rangle$. Thus $\Boot_V\subseteq \langle \iota(p)\mid p\in V \rangle$. 

The other inclusion $ \langle \iota(p)\mid p\in V \rangle \subseteq \Boot_V$ follows because $\Boot_V$ is localizing  (Lemma \ref{lemma:localizing}) and $\iota(p)\in \Boot_V$ for all $p\in V$.
\qed

\begin{lemma} \label{lemma:iota_kappa}
$\langle \iota(p) \rangle = \langle \kappa(p) \rangle$ for every $p\in \spec \Z$.
\end{lemma}
\proof
Since $\iota(0)=\kappa(0)$, we can assume that $p$ is a prime number.
There is a short exact sequence $\mathbb F_p\rightarrowtail I(p)\twoheadrightarrow I(p)$, where the first map sends $1$ to $\frac{1}{p}$ and the second map is multiplication by~$p$. This can be realized as a distinguished triangle $\Sigma\iota(p)\to \kappa(p) \to \iota(p)\to \iota(p)$ in $\Boot$. Hence $\kappa(p)\in \langle \iota(p)\rangle $, which proves ``$\supseteq$''.
For the converse inclusion, recall that for every element $x\in I(p)$ there is some positive integer $n$ such that $p^nx=0$.
This defines a filtration $0=M^0 \subseteq M^1 \subseteq \cdots \subseteq M^n \subseteq \cdots I(p)$, where $M^n$ consists of the $p^n$-torsion elements of $I(p)$. For every $n$, build an exact triangle $B^n \to B^{n+1}\to C^n \to \Sigma B^n$ with $K_*(B^{n}\to B^{n+1})=(M^n\hookrightarrow M^{n+1})$. We see from the long exact sequence in K-theory that each $K_*C^n\simeq M^{n+1}/M^n[0]$ is a graded $\mathbb F_p$-vector space. Hence $C^n\in \langle \kappa(p)\rangle$ by Cor.~\ref{cor:isotype}.
Since $B^1\simeq C^0$,  it follows by induction that $B^n\in \langle \kappa(p) \rangle$ for every $n\geq0$. Therefore, the homotopy colimit $\tilde B:= \hocolim( B^0\to B^1 \to B^2 \to \ldots ) $ (see \cite{neemanTr}*{$\S$1.6}) is also in $\langle \kappa(p)\rangle$. Since the K-theory functor $K_*$ is homological and commutes with coproducts, there is an isomorphism $I(p) = \colim_n M^n \simeq K_*\tilde B$ (\cite{neemanLoc}*{Lemma 1.5} \cite{meyernest_bc}*{Lemma 2.4}), showing that $\iota(p)\simeq \tilde B\in \langle \kappa(p)\rangle$, as required.
\qed

Therefore

\begin{prop} \label{prop:generation}
$\langle \iota(p)\mid p\in S \rangle = \langle \kappa(p)\mid p\in S \rangle$
for any subset $S\subseteq \spec \Z$.
\qed
\end{prop}

In particular, with Lemma \ref{lemma:easy_generation}, this has the following implication.

\begin{cor} \label{cor:generation_V}
$\Boot_V= \langle \kappa(p)\mid p\in V \rangle$ for any specialization closed $V\subseteq \spec \Z$.
In particular, the Bootstrap category is generated by the residue field objects:
$\Boot = \langle \kappa(p) \mid p\in \Spec \Z \rangle$. \qed
\end{cor}

\begin{remark}
Actually, the only part of the latter results that we shall need for Theorem~\ref{thm:main} is the inclusion $\C\in \langle \kappa(p) \mid p \in\Spec \Z\rangle $, which can also be established more directly by realizing in $\Boot$ the short exact sequence 
$\Z[0] \rightarrowtail \Q[0]  \twoheadrightarrow \bigoplus_p I(p)[0]$. However, this discussion allows the interested reader to compare our situation with the support theory of~\cite{bik}, cf.~$\S$\ref{section:context} below. With a little more effort, one can compute that the functors $\Gamma_p$ of \emph{loc.~cit.} correspond in $\Boot$ to $\Gamma_0=\kappa(0)\otimes-$ and $\Gamma_p = \Sigma \iota(p)\otimes-$ ($p\neq 0$).
\end{remark}

\section{Proof of Theorem~\ref{thm:main}}

\subsection{Localizing subcategories of $\Boot$}
\label{subsec:loc_subcats}
Let us prove the bijection of Theorem~\ref{thm:main}.
Thus we have to show that both compositions of the maps $S\mapsto \mathcal L_S$ and $\mathcal L\mapsto \Supp (\mathcal L)$ in the theorem yield the identity. 
In fact, it is easier to first establish the bijection where one redefines $\mathcal L_S$ the be $\langle \kappa(p)\mid p\in S \rangle$. The indentity
\begin{equation} \label{eqq}
\langle \kappa(p)\mid p\in S \rangle = \bigcap_{q\in \Spec \Z\smallsetminus S}\Ker K_*(- \otimes \kappa(q))
\end{equation}
will then easily follow: The subcategory of $\Boot$ on the right hand side is clearly localizing, so by the classification it must be generated by a unique set of residue field objects, say $\kappa(p)$ for $p\in S'$. But it is immediate from Corollary~\ref{cor:isotype} and Corollary~\ref{cor:tensor_dichotomy} that we must have $S=S'$, and therefore~\eqref{eqq}.

We shall also need the following general fact (that we have already used a couple of times), applied to $\tcat=\Boot$:

\begin{lemma} \label{lemma:ideals}
Let $\tcat$ be a triangulated category equipped with a tensor product $\otimes$ with unit object $\unit$ and which is exact and preserves coproducts. If $\tcat=\langle \unit\rangle$, then every localizing subcategory $\mathcal L\subseteq \tcat$ is a $\otimes$-ideal, i.e.: $A\otimes B\in \mathcal L$ for all $A\in \mathcal L$ and arbitrary $B\in \tcat$.
\end{lemma}

\proof
Let $\mathcal L$ be a localizing subcategory of $\tcat$. For any $A\in \mathcal L$, it is immediately checked that the subcategory $\mathcal L_A:=\{B\in \tcat \mid A\otimes B\in \mathcal L \}$ is localizing and contains~$\unit$. Hence $\tcat \subseteq \mathcal L_A$ for all $A\in \mathcal L$, as claimed. 
\qed

Now, recall that $\C \in \langle \kappa(p) \mid p\in \Spec \Z\rangle$ by Corollary~\ref{cor:generation_V}.
Since $\C$ is the tensor unit, together with Lemma~\ref{lemma:ideals} this implies 
\begin{equation} \label{hint}
\langle A \rangle= \langle A\otimes \C \rangle = \langle A\otimes \kappa(p) \mid p\in \Spec \Z\rangle
\end{equation}
 for every separable $\Cstar$-algebra~$A$. If moreover $A\in \Boot$, then
\begin{equation} \label{hint2}
A\otimes \kappa(p)\neq0
\quad 
\Rightarrow
\quad
\langle A \otimes \kappa(p) \rangle= \langle \kappa(p) \rangle 
\end{equation} 
 by Corollary~\ref{cor:residual_decomp}.
Therefore, for any localizing subcategory $\mathcal L\subseteq \Boot$ we obtain
\begin{eqnarray*}
\mathcal L \,= \, \langle A\mid A\in \mathcal L\rangle
&\stackrel{\eqref{hint}}{=}& \langle A\otimes \kappa(p) \mid A\in \mathcal L, p\in \Spec \Z\rangle \\
&\stackrel{\eqref{hint2}}{=}& \langle \kappa(p) \mid  p\in \Spec \Z  \textrm{ s.t. } \exists\, A\in \mathcal L \textrm{ with } A\otimes \kappa(p)\neq 0 \rangle \\
&=&  \mathcal L_{\Supp \mathcal L}.
\end{eqnarray*}
For the other composition, note that the inclusion $S\subseteq \Supp \mathcal L_S$ follows immediately from the definitions. On the other hand, if $q\in \spec \Z \smallsetminus S$ then $\Ker(-\otimes \kappa(q))$ is a localizing subcategory of $\Boot$ containing $\kappa(p)$ for all $p\in S$ (Cor.~\ref{cor:tensor_dichotomy}), and therefore $\Supp \mathcal L_S\subseteq S$ as well. Thus $\Supp \mathcal L_S=S$.

Thus the maps in the theorem are mutually inverse.

\subsection{Smashing subcategories of $\Boot$} 
Let us check that the bijection of Thm.~\ref{thm:main} restricts to smashing subcategories and specialization closed subsets.

\begin{lemma} \label{lemma:dichotomy_kappa}
For $\mathcal L\in \mathcal Loc(\Boot)$ and $p\in \spec \Z$, either $\kappa(p)\in \mathcal L$ or $\kappa(p)\in \mathcal L^{\perp}$.
\end{lemma}

\proof
By the classification just proved, we must have 
$\mathcal L=\mathcal L_{S}
= \langle \kappa(q)\mid q\in S \rangle$
 for some subset $S\subseteq \spec \Z$. If $\kappa(p)\not\in \mathcal L$, then $p\not\in S$. Thus $\Ext^{*}_{\Z}(\mathbb F_q,\mathbb F_p)=0$ for all $q\in S$ and therefore $\KK_*(\kappa(q),\kappa(p))=0$ by the UCT.
\qed

\begin{lemma} \label{lemma:dichotomy_iota}
If moreover $\mathcal L$ is smashing, then either $\iota(p)\in \mathcal L$ or $\iota(p)\in \mathcal L^{\perp}$.
\end{lemma}

\proof
We know from \ref{lemma:iota_kappa} that $\iota(p)\in \langle \kappa(p) \rangle$. Since $\mathcal L$ is smashing both $\mathcal L$ and its right orthogonal $\mathcal L^{\perp}$ are localizing, so the claim follows from the dichotomy of Lemma \ref{lemma:dichotomy_kappa}.
\qed

\begin{lemma} \label{lemma:smashing_p}
If $\mathcal L\subseteq \Boot$ is smashing and $\iota(0)\in  \mathcal L$, then $\iota(p)\in \mathcal L$ for all $p$. 
\end{lemma}

\proof
There is a nonzero morphism $\iota(0)\to \iota(p)$ corresponding to the surjection $\Q\to I(p)$, so $\iota(p)\not\in \mathcal L^{\perp}$. Therefore $\iota(p)\in \mathcal L$ by Lemma \ref{lemma:dichotomy_iota}.
\qed

\begin{lemma} \label{lemma:smashing}
Let $\mathcal L$ be a smashing subcategory of $\Boot$. If $0\in \Supp \mathcal L$, then $\mathcal L=\Boot$.
\end{lemma}

\proof
If $0\in \Supp \mathcal L$, then $\iota(0)\in \mathcal L$ by the first part of Thm.~\ref{thm:main}.
Since $\mathcal L$ is smashing, by Lemma \ref{lemma:smashing_p} it contains all $\iota(p)$'s and so it must be the whole $\Boot$ by Lemma \ref{lemma:easy_generation}.
\qed

In other words, the latter lemma says that if $\mathcal L$ is smashing then $\Supp \mathcal L$ is specialization closed. On the other hand, 
if $V\subseteq \spec \Z$ is specialization closed then $\mathcal L_V$ is either $\mathcal L_{\spec \Z}=\langle \C\rangle$, or $\mathcal L_V=\langle \kappa(p)\mid p\in V\rangle$ when $0\not\in V$; in both cases, it is generated by compact objects and therefore is smashing, as one can verify immediately. 

This concludes the proof.

\section{The greater context}
\label{section:context}

There are other similar classifications of localizing subcategories in the literature. Most notably, in \cite{neemanChr} A.\ Neeman classifies all localizing subcategories of the (full) derived category $\D(R)$ of a commutative noetherian ring $R$ in terms of subsets of $\Spec R$; in \cite{bik2}, Benson, Iyengar and Krause classify the localizing subcategories of the (full) stable module category $\mathsf{StMod}(kG)$ -- for $G$ any finite group and $k$ a field whose characteristic divides the order of $G$ --  in terms of subsets of the projective scheme $\mathrm{Proj}(H^*(G;k))$.

These classifications and a few others all fit into the general abstract context developed in \cite{bik, bik3}, where one works with a compactly generated triangulated category on which acts a noetherian graded commutative ring $R$. For us, $R$ would be the endomorphism ring $\Z \simeq \End_{\Boot}(\C)$ of the tensor unit object $\C\in \Boot$, but unfortunately the Bootstrap category is not compactly generated because it does not have arbitrary small coproducts. 
This is not a trivial obstruction, since not all triangular techniques employed in \emph{loc.\ cit.} can be adapted to the more general (but definitely weaker) setting of a ``compactly${}_{\aleph_1}$ generated category'' \cite{kkGarticle}*{$\S$2}, into which $\Boot$ fits. In particular, Brown representability for homological functors does not hold in $\Boot$ (see \cite{kkGarticle}*{Ex.\ 2.11}), and this is used by Bensor-Iyengar-Krause at a few crucial steps, such as their very construction of their support-defining functors $\Gamma_p$.  We don't know how much of their theory generalizes to the $\aleph$-relative situation.

Certainly however, the classification that we have just proved \emph{morally} belongs to the above-mentioned abstract context, if not factually.  Indeed, our proof is an adaptation of Neeman's original proof for the easy special case of $\D(\Z)$, as the knowledgeable reader has surely noticed.

\begin{ack}
We would like to thank the anonymous referee for suggesting a simplification in the proof which allowed us to dispense with a needlessly complicated argument involving Bousfield localizations.
\end{ack}

%


\end{document}